\documentclass[runningheads,a4paper]{llncs}

\usepackage{amssymb}
\usepackage{graphicx}

\usepackage{url}
\urldef{\mailsa}\path|{alfred.hofmann, ursula.barth, ingrid.haas, frank.holzwarth,|
\urldef{\mailsb}\path|anna.kramer, leonie.kunz, christine.reiss, nicole.sator,|
\urldef{\mailsc}\path|erika.siebert-cole, peter.strasser, lncs}@springer.com|    

\newcommand{\DS}{\LARGE \renewcommand{\baselinestretch}{1.67} \normalsize}

\usepackage{graphics} 
\usepackage{epsfig} 

\title{ Continuity-Forcing for Derivatives in Data Reconstruction }   

\author{ 
L. M. Chen}    
\institute{Department of Computer Science and Information Technology\\
University of the District of Columbia\\
lchen@udc.edu 
}

\begin{document}
\DS

\maketitle

\begin{abstract}   
The smooth function reconstruction needs to use derivatives. In 2010, we
used the gradually varied derivatives to successfully constructed
smooth surfaces for real data. We also briefly explained why 
the gradually varied derivatives are needed.  
In the this paper, we present more reasons to enlighten that 
forcing derivatives to be continuous is necessary.
This requirement seems not a must in theory for functions in continuous space, 
but it is truly important
in function reconstruction for real problems. 
This paper is also to extend the meaning of the methodology for 
gradually varied derivatives to general purposes by considering
forcing calculated derivatives to be ``continuous'' or gradually
varied.  
\end{abstract}

\section{Introduction}

In smooth data reconstruction, derivative computation is a essential 
task. However, the accuracy of derivative calculation is often
omitted in most of applications. In this paper, 
we give some examples to show the importance
of obtaining continuous derivatives. And we will discuss the 
methods to make or force the calculated derivatives to be continuous.

For a given sample point, the derivatives only affect its neighbor points
not itself. When two sample points are relatively apart from each other, 
the derivatives up to $m$-order will not change the value at the sample 
points. Therefore, we
can make continuous derivatives at the sample points to any degree
without changing the value at the sample points.

We conclude that we must force a ''continuous'' derivatives before it is used.
And such continuous must based on the definition of the space we are dealing 
with not just assume the continuous one in $R^n$. This is because
that there are always continuous interpolation from finite samples to
Euclidean space. We can use gradually varied 
function~\cite{Che89,Che92,Che10}, 
piecewise linear, or the special type of Lipschitz functions 
to restrict the meaning of continuity.



\section{Numerical Derivative Calculation and Finite Difference Methods}

In numerical mathematics, the most common method for derivatives is the 
finite difference method. The other methods some times use the idea
of finite differences. 

Let $\Delta x$ be a small constant, then

\begin{equation}
      \Delta f=f(x + \Delta x) - f(x) 
\end{equation}

\noindent We have

\begin{equation}
 f'(x) \approx \frac{\Delta f}{\Delta x} = \frac {f(x + \Delta x) - f(x)} 
                                           {\Delta x}
\end{equation}
\noindent Furthermore, 
\begin{equation}
  f''(x) \approx \frac{\Delta^2 f }{\Delta x^2}  
\end{equation}
\noindent where:

\begin{equation}
   \Delta^2 f = \Delta (\Delta f)={f(x+\Delta x) - 2 f(x) + f(x-\Delta x)}  
\end{equation}

The finite difference method may generate large error in 
calculating derivatives depends on the resolution of decomposition
of the function or space. 

We choose the following example to see 
the performance of a sampled function in discrete space:
$f(i)=(-1)^i$. This function is a sequence of $f(0)=1, f(1)= -1, 1, -1 ..., 1,
-1,...$. The first order derivative of $f$ is $f'(0)=-2, 2, ...,-2,2,...$. 
Thus, $f^{(i)} = 2^{i}$ . 

This is the worst case. We also can 
let $f(0)=L, f(1)=-L, L, -L ..., L,$
Then we will have $f= 2^{i}*L $ . 
(If we consider this function is Lipschitz~\cite{Hei05}, 
the Lipschitz constant is $2L$.)

Does the original function just have such an $m$-order derivatives?
It also depends on the sampling ratio beyond the original function. 
For instance, $Cos(x)$ is such a function if we make samples like $Cos(i*\pi)$. 
However,
$(Cos(x))'\le 1$ for all $x$. If we use the difference formula above 
to do a calculation, it must be wrong.

The better sampling for $Cos(x)$ is 1,0,-1, 0,1, 0,-1, ..., and 
$1,a, 0, -a, -1, -a, 0$, $a, 1, a, 0,-a, -1, -a ...,$ where $a=Cos(\pi/4)$.
One might say that we can do fine sampling, and there is a sampling (ratio) 
theorem. However, in general application, we do not know the periodical cycles.


\section{Methods for Continuity-Forcing of Derivatives} 

In this section, we propose to use three methods for making continuous if the derivatives are not ``continuous.'' 
In discrete space, continuity means gradually variation~\cite{Che10} or
local Lipschitz with limited Lipschitz constant locally.

\subsection{Method A: Gradually Varied Derivatives}
Let $J$ be a subset of domain $D$ and $f$ be the function from
$J$ to $\{A_1,...,A_n\}$ where $A_1<...<A_n$.   $GVF(f)$ is 
a gradually varied extension on $D$. 
$GVF(f)$ could have two meanings: (1) $f$ is already on D, ($J=D$), 
GVF(f) is a gradually 
varied uniform approximation of $f$. (2) $GVF(f)$ is gradually 
varied interpolation on $D$.  

Let $g=F^{(0)}=GVF(f)$. $g'\approx \frac {\Delta g}{\Delta x}$. Since
$\Delta x =1$, so $g'\approx {\Delta g}$. We can just use $g'= (GVF(f))'$
for simplification. Thus, $F'=(g)'=(F^{(0)})'$, so 

\[F^{(1)}=GVF(g')=GVF((F^{(0)})')=GVF(GVF'(f))\] 

Define 

$GVF^{(k)}(F^{(0)}) = GVF((GVF^{k-1}(F^{0}))')=GVF(GVF^{k-1}(F^{0}))'$. 

$GVF^{(k)}(F^{(0)})$ is an approximation of $F^{(k)}$, the $k$th-order 
of derivatives that is also continuous (the $k$th-order continuous derivatives).

In calculus, the Taylor expansion theorem states: 
A differentiable function $f$ around a point $x_0$ can be represented by a polynomial 
composed by the derivatives at the given point $x_0$. The one variable formula is as
follows: 

\begin{equation} 
\begin{array}{ll}
 f(x) = & f(x_0) + f'(x_0)(x-x_0) + \frac{f''(x_0)}{2!}(x-x_0)^2 + \cdots + \\ 
        &  \frac{f^{(k)}(x_0)}{k!}(x-x_0)^k + R_k(x)(x-x_0)^k,
\end{array}
\end{equation}

\noindent where $\lim_{x\to x_0 }R_k(x)=0$. The polynomial in (5) is 
called the Taylor polynomial or Taylor series. $R_k(x)$ is called the 
residual. 
This theorem provides us with a theoretical foundation of 
finding the $k$-th order  
derivatives at point $x_0$ since we can restore the functions around $x_0$.

After the different derivatives are obtained, we can use Taylor expansion to
update the value of the gradually varied fitted function 
(at $C^{0}$). In fact, in any order $C^{k}$, 
we can update the function using a higher order of derivatives as discussed 
in the above section.  For an $m$-dimensional space, the Taylor expansion has the 
following generalized form expanding at the point $(a_1,\dots,a_k)$:

\begin{equation}
 f(x_1,\dots,x_k) =  \sum_{n_1=0}^\infty \cdots \sum_{n_k = 0}^\infty 
                    \frac{(x_1-a_1)^{n_1}\cdots (x_k-a_k)^{n_k}}{n_1!\cdots n_k!}\,\left(\frac{\partial^{n_1 + \cdots + n_k}f}
                    {\partial x_1^{n_1}\cdots \partial x_k^{n_k}}\right)(a_1,\dots,a_k).
\end{equation}



\noindent For a function with two variables, $x$ and $y$, 
the Taylor series of the second order at expanding point $(x_0, y_0)$ is:


\begin{equation} 
\begin{array}{ll}
 f(x,y)  \approx & f(x_0,y_0) +(x-x_0)\cdot f_x (x_0,y_0) +(y-y_0)\cdot f_y(x_0,y_0) +  \\
                 &  \frac{1}{2!}[ (x-x_0)^2\cdot f_{xx}(x_0,y_0) + 2(x-x_0)(y-y_0)\cdot f_{xy}(x_0,y_0) + \\
                 &     (y-y_0)^2\cdot f_{yy}(x_0,y_0),]
\end{array} 
\end{equation}

\noindent There are several ways to implement this formula. 
For simplicity, we use $G(x)$ for $GVF(x)$. For smooth gradually varied
surface applications,  $f(x_0,y_0)$,  $f_x$, and $f_y$  are 
$G(f)$, $G^{(1)}(f)$, etc.

\begin{equation} 
\begin{array}{ll}

 f(x,y)  \approx & G(x_0,y_0) +(x-x_0)\cdot G_{x}(x_0,y_0) +(y-y_0)\cdot G_y(x_0,y_0) +  \\
                 &  \frac{1}{2!}[ (x-x_0)^2\cdot G_{xx}(x_0,y_0) + 2(x-x_0)(y-y_0)\cdot G_{xy}(x_0,y_0) + \\
                 &     (y-y_0)^2\cdot G_{yy}(x_0,y_0),]
\end{array} 
\end{equation}

The above formula shows the principle of digital-discrete 
reconstruction~\cite{Che10}.

\subsection{Method B: Derivatives Calculation by Normals on Meshed Manifolds}

Another common method for derivative calculation is to use normals.
This method is particularly used in computer graphics in computing
derivatives on meshes. After the derivatives are calculated on 
sample points, we will make the piecewise linear interpolation
on the sample points. So we get a piecewise linear derivatives. We
call this method the triangulated Derivatives.

The method for vertex (normalized) normal calculation for triangulated meshes. 
This method will do an average of the normals at each triangle around the given 
point.  Then, the final value will be normalized as the normal of the point.

We use the normal to get partial derivatives at (sample) points.
For instance, for a surface $z=f(x,y)$, the normal vector at a point $(x_0,y_0)$    
is ${\bf n}=(f_x (x_0,y_0), f_y(x_0,y_0), -1)$.    
After that, the most important step is to fit the derivatives to be gradually 
varied or ``continuous.''  
A gradually varied derivatives are necessary for further use.   
Lastly, we can have up to $k$ number of (directional) derivatives and then 
use the Taylor expansion to 
get the value surrounding the sample points. 
Again, this idea is the same as Chen's paper in 2010 ~\cite{Che10a}.
 
In general, given $m$ sample points in 3D (or higher dimension), we
want to get a smooth surface passing those points. We can first
get the piece-wise linear surface by using triangulation. So we will 
have a continue interpolation on those points.  

Given a three-dimensional surface, for a implicit form $F(x,y,z)=0$,
${\bf n}=(\nabla F)/(\sqrt(F_x^2+F_y^2+F_z^2))$. We will have $m$ normals
${\bf n_i}$  for $m$ sample points. Then use the m vectors in 3D (or higher)
 spaces to do another triangulation, we will have piece-wise linear function
 on ${\bf n_i}$. This function is the source of the first order partial derivatives
 of $F$. 

Another method for the normal and directional derivatives on surfaces is to 
 use 6 nearest points to fit a quadratic equation for 
local surface, $f(x,y)= a_0+a_1 x  + a_2 y + a_3 xy+a_4 x^2+ a_5 y^2 $,  then
we can get the normal and directional derivatives for triangles, 
any other shapes of meshes, or mesh-free cases.  
Again, we shall calculate gradually varied derivatives for this. 
Use the same method above,
so we can get the derivatives
for Taylor expansions. 

It is also true that we can use the gradually varied function on manifolds
to fit the derivative functions when the derivatives on sample points
are calculated~\cite{Che10a,CL}.  

we have used normals as the concept in 1990 to explain how $\lambda$-connected
method works.

\subsection{Method C: Lipschitz Extension of Derivatives}

This suggested method is to use Lipschitz extension to get a Lipschitz function.
At the sample points, we calculate derivatives then use
the Lipschitz extension on the derivatives at sample data location to 
get a Lipschitz derivative functions. Next to use Taylor expansions to 
obtain the whole function extension. This method is very similar to 
Method A.

For a Lipschitz function, Lipschitz continuous is not an obvious reasonable 
method for 
construction if the Lipschitz constant does not have a restriction.
the way is to find a Lipschitz construction then the 
higher order of derivatives will have less (or equal at most) Lipschitz constant. 
$Lip^{i}\le Lip^{i-1}/2$ with the same decomposition (resolution in 
grid or triangulation) is very reasonable selection.

The key to the above three methods is to always use first order derivatives 
with the adjustment in its own function. We do not use the 
second order or higher difference equations. The advantage is not
to transfer the errors of difference to higher order. In other words,
difference method is not always trustable in higher order.

\section{Discussion: The Difference From Other Techniques}

Finding continuous derivatives has been investigated by  Möller and others 
~\cite{Mol97,Mol98}. 
Their technique is to preserve calculated 
derivatives to be continuous. Our method is to force the derivatives
 to be continuous
when the derivatives are not continuous. 

These two methods have philosophical difference in theory and practice. 
Moller's work is a specific technique. One can ask:
 how to obtain a general method for most of cases when
the function formula is not known?

Our method is a general method for most of cases when sample points are 
provided.  We construct the derivative points in the function,
then make them to be continuous. That is the difference.

On the other hand,
a popular method called subdivision method for smooth function reconstruction 
in graphics~\cite{CC}.  The final function is $C^{(1)}$ or  $C^{(2)}$. Since this method is not designed
for data reconstruction, the shape of the final function must be designed. Without dense
sample points, the result of the method is unknown. However, our method presented
in the above sections has fitted all reasonable data points on every point in the manifold.
Those points are now perfect as the sample points for the second level of fitting--
using subdivision methods. We are working on the implementation of this part. We
believed that we would obtain very good results by keeping the original sample points
unchanged (the method in 2D is sometimes called 4 point interpolation). The key is to 
add some points around the corner or ridge to make the surface repeatedly smoother . 

Moving least square methods can also be applied. Moving least square is a 
mesh-free method~\cite{LS}. There is no need to partition the domain. It uses a normal 
distribution like function to weight the sample data for a local polynomial fitting.
The fitting is followed by the principle of the least square. While we move from
one point to another, the method will fit a new value to the current point. 
This method must rely on the dense sample points with a balanced distribution.
There may be no points near the fitting point. In another
situation, there may be many points around a fitting point. How do we select
a weight equation?  
The disadvantage of this method is that we need some knowledge for the weight equation, 
the lengths of circles for the weight changes. An artificial intelligence method may be
needed for this determination.  The process of our harmonic data reconstruction will
provide the first good fitted data. If we want to use a local
polynomial to fit the data again, we can use moving least square methods at the 
top of our algorithm. 

Smooth gradually varied functions use gradually varied derivatives and the 
Taylor expansion for the function reconstruction. Gradually varied derivatives
are needed because finite difference might not get continuous derivatives. Using  
gradually varied fitting on the derivatives is necessary in such a case.
In~\cite{Che10a}, we presented this method for a 2D rectangle domain. For 
2D manifolds, the calculation of derivatives cannot just use finite difference
methods.   We can use Method B to solve the problem.

\end{document}